# Numerical modeling of wing crack propagation accounting for fracture contact mechanics

Hau Dang-Trung[*], Eirik Keilegavlen, Inga Berre

*Department of Mathematics, University of Bergen, Bergen, Norway*

ABSTRACT

As a consequence of shearing, wing cracks can emerge from pre-existing fractures. The process involves the interaction of sliding of the existing fracture surfaces and the tensile material failure that creates wing cracks. This work devises a numerical model to investigate how wing cracks emerge, propagate and connect pre-existing fractures under shear processes. A mathematical and numerical model for wing crack propagation based on linear elastic fracture mechanics that also accounts for fracture contact mechanics is presented. Computational efficiency is ensured by an adaptive remeshing technique. The numerical model is verified and validated through a comparison of the analytical and experimental results. Additional numerical examples illustrate the performance of the method for complex test cases where wing-cracks develop for multiple pre-existing and interacting fractures.

[*]Corresponding author: Hau Dang-Trung, Department of Mathematics, University of Bergen, Bergen, Norway.

Email: Hau.Dang@uib.no; dtrhau@gmail.com

## 1. Introduction

Wing cracks can develop from a pre-existing fracture when the fracture is subjected to shear processes. This occurs for many applications where fractured media are subjected to anisotropic stress regimes. For example, in fractured subsurface systems, fractures will slip if shear forces overcome the cohesion and frictional strength of the contact between the fracture surfaces. This can occur due to natural changes in tectonic stresses, but the process can also be induced by fluid injection, such as in situations of geothermal reservoirs. In the latter case, elevated pressures reduce the effective normal stress on the fracture, ultimately causing slip if the reduction in the normal stress is sufficient for the shear forces to overcome the cohesion and frictional resistance of the fracture. The slip of the fracture surfaces in opposite directions can cause the fracture to propagate in the form of wing cracks, possibly creating enhanced reservoir connectivity (Cheng et al., 2019; Jung, 2013; McClure and Horne, 2014; Norbeck et al., 2018). Understanding this mechanism is, thus, crucial in the simulation of fractured subsurface formations.

Many experimental studies have been published that consider the formation, growth and connection of wing cracks caused by external compressive loading in specimens made of rock or rock-like materials (Haeri et al., 2014a, 2014b; Horii and Nemat-Nasser, 1985; Ingraffea and Manu, 1980). In these experiments, if the pre-existing fracture is not perpendicular to the external load, wing cracks emerge at the tip and tend to align with the direction of the maximum compressive stress. The same conclusion is drawn from mathematical modeling. Based on the finite element method (FEM), Ingraffea and Heuze (1980) predicted the propagation of wing cracks in rock structures by using three different criteria based on stress, energy and strain. Primary crack trajectories predicted by the stress and energy criteria are in good agreement with the observed trajectories. Based on the phase-field model (Bryant and Sun, 2018) and a modified phase-field model (Zhang et al., 2017), wing crack propagation was modeled using energy criteria that divided the active energy density into distinct parts corresponding to different crack modes (mode I and mode II). Sharafisafa and Nazem (2014) used the vector level set with both the discrete element method (DEM) and the extended finite element method (XFEM) to model the wing crack propagation and coalescence in fractured rock masses. Among these methods and failure criteria, the FEM is the simplest method in implementation and the stress criterion is one of the most extensively used and least complicated. Specifically, the combination of FEM with the stress criterion (the maximum tangential stress) for modeling of wing cracks has been shown to produce simulation in good agreement with observed crack trajectories (Gonçalves da Silva and Einstein, 2013; Ingraffea and Heuze, 1980).

While wing cracks develop as tensile fractures, the pre-existing fractures that wing cracks emerge from may be either open or in contact. This necessitates the inclusion of fracture contact mechanics in the wing crack models (Hüeber et al., 2008; Oden and Pires, 1983). There are several ways to formulate the contact mechanics corresponding to the different types of discretizations. For example, Kim and Duarte (2015) simulated the mode I propagation of cohesive fractures in 3D by the cohesive law using a generalized finite element method. Because of the difference of material behaviour in the vicinity of propagating crack fronts compared to the rest of the domain, this approach requires updates of global-local enrichments during the analysis, which increases the computational cost. Hesch et al. (2016) formulated the contact mechanics by Coulomb's friction law and Karush Kuhn-Tucker (KKT) conditions applied to a phase-field approach within the context of isogeometric analysis. A fourth order approach for the crack-density functional was used to ensure sufficient accuracy of the chosen phase-field approach. This leads to a request of at least $C^1$ continuity within the domain. Also using Coulomb's friction law and the KKT conditions, Nejati et al. (2016) modeled the internal contact in fractured media by a sophisticated algorithm based on isoparametric integration-point-to-integration-point discretization of the contact contribution to enforce the contact constraint accurately over the crack surfaces. Based on the semianalytical displacement discontinuity method, Kamali and Ghassemi (2018) developed a simulation model in which the closed natural fractures were represented by so-called contact displacement discontinuity elements (Asgian, 1988), approximating the contact mechanics condition. However, the approach has limitations in dealing with the interaction between multiple fractures due to inherent limitations of the semianalytical displacement discontinuity method.

The FEM model for elasticity can be derived by using one of the following widely used methods: the weighted-residual method based on the linear momentum balance or minimization of an energy function (Liu and Quek, 2003). By using the energy principle, the contact should be considered as an inequality constraint of the optimization formulation of the potential energy. This means that the potential energy is minimized while satisfying a contact constraint assumed to be a nonpenetration condition between the surfaces of the fracture. The inequality constraint can be solved by some methods, such as the active-set, Frank-Wolfe, penalty or barrier methods (Hüeber et al., 2008; Hüeber and Wohlmuth, 2005).

An inherent problem in the simulation of fracture propagation is the disparate length scales. While the simulation domain can be quite large, the fracturing processes occur on a scale that is several orders of magnitude smaller. Moreover, most numerical methods for fracture propagation are dependent on resolving the fracture in a grid; however, the fracture path is not known a priori. A possible remedy for both of these issues is to apply adaptive remeshing (ARM) techniques to refine and adjust the mesh around an advancing fracture path.

This paper presents a mathematical model and corresponding numerical solution approach to simulate the development of wing cracks while accounting for fracture contact mechanics. First, in Section 2, the mathematical model for wing crack propagation is formulated based on the linear elasticity theory, in combination with the criteria for a mixed-mode fracture propagation. Fracture surfaces are allowed to be in contact or fully open, modeled by contact mechanics formed by the KKT conditions. Section 3 presents the numerical solution approach. The governing equations are discretized using a finite element method with collapsed quarter-point elements at the fracture tips. This is combined with an adaptive remeshing technique based on error estimates and Laplacian smoothing. The contact mechanics are implemented by using an active set method. Section 4 presents several numerical test cases. The obtained results are compared with both the analytical and experimental data to verify, validate and show the accuracy of the proposed model and procedure. Finally, more complex test cases where wing cracks develop for multiple pre-existing and interacting fractures show the capability of the proposed approach in modeling the development of wing cracks under shear processes.

**Acknowledgements**

This work was funded by the ERiS project, grant #267909, Research Council of Norway.

Please cite this article as: H. Dang-Trung, E. Keilegavlen and I. Berre, 2020. Numerical modeling of wing crack propagation accounting for fracture contact mechanics. International Journal of Solids and Structures xx, xx-xx
DOI: 10.1016/j.ijsolstr.2020.08.017

## References

Anderson, T.L., 2017. Fracture Mechanics: Fundamentals and Applications, Fourth Edition. CRC Press.

Asgian, M.I., 1988. A numerical study of fluid flow in a deformable, naturally fractured reservoir: The influence of pumping rate on reservoir response, in: The 29th U.S. Symposium on Rock Mechanics (USRMS). American Rock Mechanics Association, Minneapolis, Minnesota, p. 8.

Atkinson, C., Smelser, R.E., Sanchez, J., 1982. Combined mode fracture via the cracked Brazilian disk test. Int. J. Fract. 18, 279–291.

Barsoum, R.S., 1977. Triangular quarter-point elements as elastic and perfectly-plastic crack tip elements. Int. J. Numer. Methods Eng. 11, 85–98.

Barsoum, R.S., 1976. On the use of isoparametric finite elements in linear fracture mechanics. Int. J. Numer. Methods Eng. 10, 25–37.

Bobet, A., Einstein, H.H., 1998. Fracture coalescence in rock-type materials under uniaxial and biaxial compression. Int. J. Rock Mech. Min. Sci. 35, 863–888.

Bryant, E.C., Sun, W., 2018. A mixed-mode phase field fracture model in anisotropic rocks with consistent kinematics. Comput. Methods Appl. Mech. Eng. 342, 561–584.

Buell, W.R., Bush, B.A., 1973. Mesh Generation-A Survey. J. Eng. Ind. 95, 332–338.

Chen, L.S., Kuang, J.H., 1992. A modified linear extrapolation formula for determination of stress intensity factors. Int. J. Fract. 54, R3–R8.

Cheng, Q., Wang, X., Ghassemi, A., 2019. Numerical simulation of reservoir stimulation with reference to the Newberry EGS. Geothermics 77, 327–343.

Erdogan, F., Sih, G.C., 1963. On the Crack Extension in Plates Under Plane Loading and Transverse Shear. J. Basic Eng. 85, 519–525.

Field, D.A., 1988. Laplacian smoothing and Delaunay triangulations. Commun. Appl. Numer. Methods 4, 709–712.

Gonçalves da Silva, B., Einstein, H., 2013. Modeling of crack initiation, propagation and coalescence in rocks. Int. J. Fract. 182, 167–186.

Haeri, H., Shahriar, K., Fatehi Marji, M., Moarefvand, P., 2014a. Experimental and numerical study of crack propagation and coalescence in pre-cracked rock-like disks. Int. J. Rock Mech. Min. Sci. 67, 20–28.

Haeri, H., Shahriar, K., FatehiMarji, M., Moarefvand, P., 2014b. On the crack propagation analysis of rock like Brazilian disc specimens containing cracks under compressive line loading. Lat. Am. J. Solids Struct. 11, 1400–1416.

Henshell, R.D., Shaw, K.G., 1975. Crack tip finite elements are unnecessary. Int. J. Numer. Methods Eng. 9, 495–507.

Hesch, C., Franke, M., Dittmann, M., Temizer, I., 2016. Hierarchical NURBS and a higher-order phase-field approach to fracture for finite-deformation contact problems. Comput. Methods Appl. Mech. Eng. 301, 242–258.

Horii, H., Nemat-Nasser, S., 1985. Compression-induced microcrack growth in brittle solids: Axial splitting and shear failure. J. Geophys. Res. Solid Earth 90, 3105–3125.

Hüeber, S., Stadler, G., Wohlmuth, B., 2008. A Primal-Dual Active Set Algorithm for Three-Dimensional Contact Problems with Coulomb Friction. SIAM J. Sci. Comput. 30, 572–596.

Hüeber, S., Wohlmuth, B., 2005. A primal–dual active set strategy for non-linear multibody contact problems. Comput. Methods Appl. Mech. Eng. 194, 3147–3166.

Ingraffea, A.R., Heuze, F.E., 1980. Finite element models for rock fracture mechanics. Int. J. Numer. Anal. Methods Geomech. 4, 25–43.

Ingraffea, A.R., Manu, C., 1980. Stress-intensity factor computation in three dimensions with quarter-point elements. Int. J. Numer. Methods Eng. 15, 1427–1445.

Jaeger, J.C., Cook, N.G.W., Zimmerman, R., 2007. Fundamentals of Rock Mechanics. Wiley.

Jung, R., 2013. EGS — Goodbye or Back to the Future, in: Effective and Sustainable Hydraulic Fracturing. InTech, pp. 95–121.

Kachanov, M., 1987. Elastic solids with many cracks: A simple method of analysis. Int. J. Solids Struct. 23, 23–43.

Kamali, A., Ghassemi, A., 2018. Analysis of injection-induced shear slip and fracture propagation in geothermal reservoir stimulation. Geothermics 76, 93–105.

Khoei, A.R., Azadi, H., Moslemi, H., 2008. Modeling of crack propagation via an automatic adaptive mesh refinement based on modified superconvergent patch recovery technique. Eng. Fract. Mech. 10, 2921–2945.

Kim, J., Duarte, C.A., 2015. A new generalized finite element method for two-scale simulations of propagating cohesive fractures in 3-D. Int. J. Numer. Methods Eng. 104, 1139–1172.

Kuang, J.H., Chen, L.S., 1993. A displacement extrapolation method for two-dimensional mixed-mode crack problems. Eng. Fract. Mech. 46, 735–741.

Laures, J.-P., Kachanov, M., 1991. Three-dimensional interactions of a crack front with arrays of penny-shaped microcracks. Int. J. Fract. 48, 255–279.

Legrand, L., Lazarus, V., 2015. Front shape and loading evolution during cracks coalescence using an incremental perturbation method. Eng. Fract. Mech. 133, 40–51.

Liu, G.R., Quek, S.S., 2003. Finite Element Method: A Practical Course. Butterworth-Heinemann, Oxford, pp. 1–11.

McClure, M.W., Horne, R.N., 2014. An investigation of stimulation mechanisms in Enhanced Geothermal Systems. Int. J. Rock Mech. Min. Sci. 72, 242–260.

Mirza, F.A., Olson, M.D., 1978. Energy convergence and evaluation of stress intensity factor $K_I$ for stress singular problems by mixed finite element method. Int. J. Fract. 14, 555–573.

Nejati, M., Paluszny, A., Zimmerman, R.W., 2016. A finite element framework for modeling internal frictional contact in three-dimensional fractured media using unstructured tetrahedral meshes. Comput. Methods Appl. Mech. Eng. 306, 123–150.

Norbeck, J.H., McClure, M.W., Horne, R.N., 2018. Field observations at the Fenton Hill enhanced geothermal system test site support mixed-mechanism stimulation. Geothermics 74, 135–149.

Oden, J.T., Pires, E.B., 1983. Nonlocal and Nonlinear Friction Laws and Variational Principles for Contact Problems in Elasticity. J. Appl. Mech. 50, 67–76.

Paluszny, A., Matthäi, S.K., 2009. Numerical modeling of discrete multi-crack growth applied to pattern formation in geological brittle media. Int. J. Solids Struct. 46, 3383–3397.

Paris, P., Erdogan, F., 1963. A Critical Analysis of Crack Propagation Laws. J. Fluids Eng. 85, 528–533.

Parks, D.M., 1974. A stiffness derivative finite element technique for determination of crack tip stress intensity factors. Int. J. Fract. 10, 487–502.

Pin, T., Pian, T.H.H., 1973. On the convergence of the finite element method for problems with singularity. Int. J. Solids Struct. 9, 313–321.

Renshaw, C.E., Pollard, D.D., 1994. Numerical simulation of fracture set formation: A fracture mechanics model consistent with experimental observations. J. Geophys. Res. Solid Earth 99, 9359–9372.

Sharafisafa, M., Nazem, M., 2014. Application of the distinct element method and the extended finite element method in modelling cracks and coalescence in brittle materials. Comput. Mater. Sci. 91, 102–121.

Tada, H., Paris, P.C., Irwin, G.R., 2000. The stress analysis of cracks handbook. ASME Press, Norwood Mass.

Thomas, R.N., Paluszny, A., Zimmerman, R.W., 2017. Quantification of Fracture Interaction Using Stress Intensity Factor Variation Maps. J. Geophys. Res. Solid Earth 122, 7698–7717.

Treifi, M., Olutunde Oyadiji, S., Tsang, D.K.L., 2008. Computations of modes I and II stress intensity factors of sharp notched plates under in-plane shear and bending loading by the fractal-like finite element method. Int. J. Solids Struct. 45, 6468–6484.

Wohlmuth, B., 2011. Variationally consistent discretization schemes and numerical algorithms for contact problems. Acta Numer. 20, 569–734.

Wong, L.N.Y., Einstein, H.H., 2009. Systematic evaluation of cracking behavior in specimens containing single flaws under uniaxial compression. Int. J. Rock Mech. Min. Sci. 46, 239–249.

Zhang, X., Sloan, S., Vignes, C., Sheng, D., 2017. A modification of the phase-field model for mixed mode crack propagation in rock-like materials. Comput. Methods Appl. Mech. Eng. 322, 123–136.

Zienkiewicz, O.C., Taylor, R.L., Zhu, J.Z., 2005. The Finite Element Method: Its Basis and Fundamentals. Elsevier Science.

Zienkiewicz, O.C., Zhu, J.Z., 1987. A simple error estimator and adaptive procedure for practical engineerng analysis. Int. J. Numer. Methods Eng. 24, 337–357.